\documentclass[10pt]{article}

\usepackage{amsmath}
\usepackage{amsthm}
\usepackage{amssymb}
\usepackage{amscd}
\usepackage{amsfonts}  
\usepackage{makeidx}
\usepackage{hyperref}

\newcommand{\bbC}{{\mathbb C}}
\newcommand{\bbQ}{{\mathbb Q}}
\newcommand{\bbZ}{{\mathbb Z}}

\newcommand{\calfin}{{\mathcal F}{\mathcal I}{\mathcal N}} 
\newcommand{\caln}{{\mathcal N}}

\newcommand{\colim}{\operatorname{colim}}
\newcommand{\Ext}{\operatorname{Ext}}
\newcommand{\id}{\operatorname{id}}
\newcommand{\res}{\operatorname{res}}

\newcommand{\Tor}{\operatorname{Tor}}

\newtheorem{theorem}{Theorem}
\newtheorem{lemma}[theorem]{Lemma}
\newtheorem{conjecture}[theorem]{Conjecture}

\newcounter{commentcounter}
\newcommand{\comment}[1]                      
{\stepcounter{commentcounter}
{\bf Comment \arabic{commentcounter}}: {\ttfamily #1} }

\begin{document}

\title{Groups of small homological dimension and the Atiyah Conjecture}

\author{Peter Kropholler, Peter Linnell  and Wolfgang L\"uck}

\maketitle


\typeout{-----------------------  Abstract  ------------------------}

\begin{abstract}

A group has homological dimension $\le 1$ if it is locally free.
We prove the converse provided that $G$ satisfies the Atiyah Conjecture about
$L^2$-Betti numbers. We also show that 
a finitely generated elementary amenable group $G$ of
cohomological dimension $\le 2$ possesses a finite
$2$-dimensional model for $BG$ and in particular that $G$ is
finitely presented and 
the trivial $\bbZ G$-module $\bbZ$ has
a $2$-dimensional resolution by finitely generated free
$\bbZ G$-modules.

\smallskip
\noindent
Key words: (co-)homological dimension, von Neumann dimension, Atiyah Conjecture.\\
Mathematics Subject Classification 2000: 18G20, 46L99.
\end{abstract}


\section{Notation}

Throughout this paper let $G$ be a (discrete) group. It has 
\emph{homological dimension $\le n$} if 
$H_p(G;M) = \Tor^{\bbZ G}_p(\bbZ,M)$ vanishes for each $\bbZ G$-module
$M$ and each $p >n$. It has 
\emph{cohomological dimension $\le n$} if 
$H^p(G;M) = \Ext_{\bbZ G}^p(\bbZ,M)$ vanishes for each $\bbZ G$-module
$M$ and each $p > n$.

We call $G$ \emph{locally free} if each finitely
generated subgroup is free. The \emph{class of elementary amenable groups}
is defined as the smallest
class of groups, which  contains all finite and all abelian groups
and is closed under taking subgroups, taking quotient groups,
extensions and directed unions. Each elementary amenable group is amenable,
but the converse is  not true.


\section{Review of the Atiyah Conjecture}

Denote by $\caln(G)$ the group von Neumann algebra associated to $G$ which we
will view as a ring (not taking the topology into account) throughout this paper.
For a $\caln(G)$-module $M$ let $\dim_{\caln(G)}(M) \in [0,\infty]$ be its dimension in
the sense of \cite[Theorem 6.7]{Lueck(2002)}. Let $\frac{1}{\calfin(G)}\bbZ \subseteq \bbQ$ be the
additive abelian subgroup of $\bbQ$ generated by the inverses $|H|^{-1}$ of the orders
$|H|$ of finite subgroups $H$ of $G$. Notice that 
$\frac{1}{\calfin(G)}\bbZ$ agrees with $\bbZ$ if and only if $G$ is torsion-free.

\begin{conjecture}[Atiyah Conjecture] \label{con:_Atiyah_Conjecture}
Consider a ring $A$ with $\bbZ \subseteq A \subseteq \bbC$. 
The \emph{Atiyah Conjecture for $A$ and $G$}
says that for each finitely presented $AG$-module $M$ we have
$$\dim_{\caln(G)}\left(\caln(G) \otimes_{AG} M\right) \in \frac{1}{\calfin(G)}\bbZ.$$
\end{conjecture} 
For a discussion of this conjecture and the classes of groups for which it is known
we refer for instance to \cite[Section 10.1]{Lueck(2002)}.
It is not clear whether the Atiyah conjecture is subgroup closed;
however in the case $G$ is torsion-free, then it certainly is.  This
can be seen from \cite[Theorem 6.29(2)]{Lueck(2002)}.
We mention Linnell's result \cite{Linnell(1993)} that the 
Atiyah Conjecture
is true for $A = \bbC$ and all groups $G$
which can be written as an extension with a free group as kernel
and an elementary amenable group as quotient and possess an upper bound
on the orders of its finite subgroups. The Atiyah Conjecture
has also been proved by Schick~\cite{Schick(2000)} for $A = \bbQ$ and torsion-free groups $G$
which are residually torsion-free elementary amenable.


\section{Results}

\begin{theorem} \label{the:_characterization_of_hdim_le_1}
A locally free group $G$ has homological dimension $\le 1$.

If $G$ is a group of homological dimension $\le 1$ and the 
Atiyah Conjecture~\ref{con:_Atiyah_Conjecture} holds for $G$, then
$G$ is locally free.
\end{theorem}

\begin{theorem} 
\label{the_finitely_generated_elementary_amenable_groups_with_cd_le_2_and-2-dim_res}
Let $G$ be an elementary amenable group of cohomological dimension $\le
2$. Then

\begin{enumerate}

\item \label{the_finitely_generated_elementary_amenable_groups_with_cd_le_2_and-2-dim_res: fin  gen.}
Suppose that $G$ is finitely generated.  Then $G$  possesses a presentation of the form
$$\langle x,y \mid yxy^{-1} = x^n \rangle.$$
In particular there is a finite $2$-dimensional model for $BG$ and the
trivial $\bbZ G$-module $\bbZ$ possesses a $2$-dimensional resolution 
by finitely generated free $\bbZ G$-modules;

\item  \label{the_finitely_generated_elementary_amenable_groups_with_cd_le_2_and-2-dim_res: not fin.gen}
Suppose that $G$ is countable but not finitely generated. Then $G$ is a non-cyclic subgroup of
the additive group $\bbQ$.

\end{enumerate}
\end{theorem}


\section{Proofs}

\begin{lemma} \label{lem:_projectives_in_finitely_generated_free_modules}
Let $A$ be a ring with $\bbZ \subset A \subset \bbC$.
Let $P$ be a projective $AG$-module
such that for some finitely generated $AG$-submodule $M \subset P$ we have
$\dim_{\caln(G)}(\caln(G)\otimes_{AG} P/M) = 0$. Then $P$ is finitely generated.
\end{lemma}
{\bf Proof:} Choose a free $AG$-module $F$ and  $AG$-maps 
$i\colon P \to F$ and $r \colon F \to P$
with $r \circ i = \id$. Since $M \subset P$ is finitely generated, there is
a finitely generated free direct summand $F_0 \subset F$  with $i(M)
\subset F_0$ and $F_1:=F/F_0$ a free $AG$-module.
Hence $i$ induces a map $f\colon P/M \to F_1$.
It suffices to show that $f$ is trivial because then $i(P) \subset F_0$ and 
the restriction of $r$ to $F_0$ yields an epimorphism $F_0 \to P$.

Let $g\colon AG \to P/M$ be any $AG$-map.  The map
$\caln(G) \otimes_{AG} (f \circ g)$ factorizes through $\caln(G) \otimes_{AG} P/M$.
Hence its image has von Neumann dimension zero because $\dim_{\caln(G)}$
is additive \cite[Theorem 6.7]{Lueck(2002)} and 
$\dim_{\caln(G)}(\caln(G)\otimes_{AG} P/M) = 0$ holds by assumption. Since the von Neumann algebra
$\caln(G)$ is semi-hereditary (see \cite[Theorem 6.5 and Theorem 6.7]{Lueck(2002)}), 
the image of $\caln(G) \otimes_{AG} (f \circ g)$
is a finitely generated projective $\caln(G)$-module, whose von Neumann dimension is zero,
and hence is the zero-module. Therefore  $\caln(G) \otimes_{AG} (f \circ g)$ is the zero map.
Since $AG \to \caln(G)$ is injective, $f \circ g$ is trivial. 
This implies that $f$ is trivial since $g$ is any $AG$-map. \qed

\begin{lemma} \label{lem:_caln(G)_otimes_P/M_has_zero_dimension}
Let $A$ be a ring with $\bbZ \subset A \subset \bbC$.
Suppose that there is a positive integer $d$ such that the order of any finite subgroup
of $G$ divides $d$ and that the Atiyah Conjecture holds for $A$ and $G$. Let $N$
be a $AG$-module. Suppose that $\dim_{\caln(G)}(\caln(G) \otimes_{AG} N) < \infty$. Then there is
a finitely generated $AG$-submodule $M \subset N$ with 
$\dim_{\caln(G)}(\caln(G) \otimes_{AG} N/M) = 0$.
\end{lemma}
{\bf Proof:} Since $N$ is the colimit of the directed system of its finitely generated $AG$-modules
$\{M_i \mid i \in I\}$ and tensor products commute with colimits, we get
$\colim_{i \in I}\caln(G) \otimes_{AG} N/M_i = 0$. Additivity (see
\cite[Theorem 6.7]{Lueck(2002)}) implies 
$\dim_{\caln(G)}(\caln(G) \otimes_{AG} N/M_i) < \infty$ for all $i \in
I$ since $\dim_{\caln(G)}(\caln(G) \otimes_{AG} N) < \infty$ holds by assumption.
We conclude from Additivity and Cofinality (see \cite[Theorem
6.7]{Lueck(2002)}) and the fact that the functor
colimit over a directed system of modules is exact
$$\inf\{\dim_{\caln(G)}(\caln(G) \otimes_{AG} N/M_i)\mid i \in I \} = 0.$$
The assumption about $G$ implies using \cite[Lemma 10.10 (4)]{Lueck(2002)}
$$d \cdot \dim_{\caln(G)}(\caln(G) \otimes_{AG} N/M_i) \in \bbZ.$$
Hence there must be an index $i \in I$ with 
$\dim_{\caln(G)}(\caln(G) \otimes_{AG} N/M_i)=0$.
\qed

\bigskip
{\bf Proof of Theorem~\ref{the:_characterization_of_hdim_le_1}:}
A finitely generated free group has obviously homological dimension $\le 1$.
Since homology is compatible with colimits over directed systems (in contrast to
cohomology), we get for every  group $G$, which is the directed union of the family of subgroups
$\{G_i \mid  i \in I\}$, and every $\bbZ G$-module $M$
$$H_n(G;M) = \colim_{i \in I} H_n(G_i;\res_i M),$$
where $\res_i M$ is the restriction of $M$ to a $\bbZ G_i$-module.
Hence any locally free group has homological dimension $\le 1$.

Suppose that $G$ has homological dimension $\le 1$.
Let $H \subset G$ be a finitely generated subgroup. 
Then the homological dimension of
$H$ is $\le 1$. 
Since each countably presented flat module is of projective dimension 
$\le 1$ \cite[Lemma 4.4]{Bieri(1981)}, we conclude
that the  cohomological dimension of $H$ is $\le 2$.
We can choose an exact sequence $0 \to P \to \bbZ H^s \to \bbZ H \to \bbZ$,
where $s$ is the number of generators and $P$ is projective. Since
the homological dimension is $\le 1$, the induced map
$\caln(H)\otimes_{\caln(H)} P \to \caln(H)\otimes_{\caln(H)} \bbZ H^s$ is injective and
hence $\dim_{\caln(H)}(\caln(H)\otimes_{\caln(H)} P) \le
\dim_{\caln(H)}(\caln(H)\otimes_{\caln(H)} \bbZ H^s) = s$. Suppose that
$G$ satisfies the Atiyah Conjecture.
Since $G$ cannot contain a non-trivial finite subgroup, $H$ also
satisfies the Atiyah Conjecture, and
Lemma~\ref{lem:_projectives_in_finitely_generated_free_modules}  and 
Lemma~\ref{lem:_caln(G)_otimes_P/M_has_zero_dimension} imply that $P$ is finitely generated.
Hence $H$ is of type $FP$.
Since each  finitely presented flat module is projective
 \cite[Lemma 4.4]{Bieri(1981)}, the cohomological dimension and 
the homological dimension agree for
groups of type $FP$.
Hence $H$ has cohomological dimension $1$.
A result of Stallings~\cite{Stallings(1968)} implies that $H$ is free. \qed
\bigskip

In \cite{Hillman91} the notion of Hirsch length for an elementary
amenable group was defined, generalizing that of the Hirsch length of
a solvable group.  This was used in the proof of
\cite[Corollary 2]{HillmanLinnell92} to show that an elementary
amenable group of finite cohomological dimension is virtually solvable
with finite Hirsch number, see \cite[Theorem 1.11]{Hillman02} for
further details.  We can now state

\begin{lemma} 
\label{lem:_EAgroups_of_hdim2}
If $G$ is an elementary amenable group of homological dimension
$\le 2$, then $G$ is metabelian.
\end{lemma}

{\bf Proof:}
A group is metabelian if and only if each finitely generated subgroup
is metabelian. Hence we can assume without loss of generality
that $G$ is finitely generated.  Then
by the above remarks and \cite[Theorem 7.10(a)]{Bieri(1981)},
$G$ is virtually solvable of Hirsch length $\le 2$.

If $G$ has Hirsch length 1, then $G$ is infinite
cyclic, so we may assume that $G$ has Hirsch length 2.
Let $N$ denote the Fitting subgroup of $G$ (so
$N$ is generated by the nilpotent normal subgroups of $G$ and is a
locally nilpotent normal subgroup).

Suppose that $N$ has finite index in $G$. Then $N$
is finitely generated and is therefore free abelian of rank 2.  Also
$G/N$ acts faithfully by conjugation on $N$ (a torsion-free group
with a central subgroup of finite index must be abelian).
If $g \in G \setminus
N$, then $g^r \in N \setminus 1$ for some positive integer $r$ and
thus $g$ fixes a nonidentity element of $N$.  We deduce that $|G/N| \le
2$ and it follows that $G$ is metabelian.  

On the other hand if $N$ has infinite index in $G$, 
then it has Hirsch length 1. Hence every finitely generated subgroup of
$N$ is trivial or isomorphic to $\bbZ$. This implies that $N$ is
abelian and any automorphism of finite order $f \colon N \to N$
has the property that $f(x) \in \{x,-x\}$ holds for $x \in N$.  Since the group
$G/N$ acts faithfully by conjugation on $N$ and is virtually cyclic, we conclude that
$G/N$ is isomorphic to $\bbZ$ or $\bbZ \times \bbZ/2$. Hence
$G$ is metabelian.
\qed

\bigskip
{\bf Proof of Theorem~\ref{the_finitely_generated_elementary_amenable_groups_with_cd_le_2_and-2-dim_res}:}
Since $G$ has cohomological dimension 2, it certainly has homological
dimension at most 2 and so by Lemma \ref{lem:_EAgroups_of_hdim2} is
metabelian.  A result of Gildenhuys \cite[Theorem 5]{Gildenhuys(1979)},
states that a solvable group $G$ of cohomological dimension 2 has a
presentation of the form $\langle x,y ; y^{-1}xy = x^n \rangle$ for
some $n\in \mathbb {Z}$ if $G$ is finitely generated and
is a non-cyclic subgroup of the additive group $\bbQ$ if $G$ is not
finitely generated. Given a  torsion free finitely generated
one-relator group $G$, the finite two-dimensional $CW$-complex
associated to a presentation with finitely many generators and one
non-trivial relation is a model for $BG$
(see \cite[Chapter III \S\S 9 -11]{Lyndon-Schupp(1977)}).
This finishes the proof of
Theorem~\ref{the_finitely_generated_elementary_amenable_groups_with_cd_le_2_and-2-dim_res}.
\qed

\typeout{-------------------------- references -----------------------}

\bigskip
{\bf Address:}\\[1mm]
Peter Kropholler, Department of Mathematics, University at Glasgow, University Garden,
Glasgow G12 8QW, Scotland, p.h.kropholler@maths.gla.ac.uk,
http://www.maths.gla.ac.uk/people/?id=289
\\[1mm]
Peter Linnell, Department of Mathematics, Virginia Tech, Blacksburg, VA 24061-0123, USA,
linnell@math.vt.edu, http://www.math.vt.edu/people/linnell/
\\[1mm]
Wolfgang L\"uck, FB Mathematik, Universit{\"a}t M{\"u}nster,
Einsteinstr. 62, D-48149 M{\"u}nster, Germany,
lueck@math.uni-muenster.de, http://wwwmath.uni-muenster.de/math/u/lueck/


\begin{thebibliography}{99}

\bibitem{Bieri(1981)} {\bf Bieri, R}:
\emph{``Homological dimension of discrete groups''}, 2-nd edition,
Queen Mary College Mathematics Notes, Mathematics Department, Queen Mary College, London
(1981).

\bibitem{Gildenhuys(1979)} {\bf Gildenhuys, D.}:
\emph{``Classification of soluble groups of cohomological dimension''},
Math. Z. 166, 21--25 (1979).

\bibitem{Hillman91}
{\bf Hillman, Jonathan~A.}:
\newblock \emph{``Elementary amenable groups and {$4$}-manifolds with {E}uler
  characteristic 0"},
\newblock {\em J. Austral. Math. Soc. Ser. A}, 50(1):160--170, 1991.

\bibitem{Hillman02}
{\bf  Hillman, J.~A.}:
\newblock {\em ``Four-manifolds, geometries and knots"}, volume~5 of {\em Geometry
  \& Topology Monographs}.
\newblock Geometry \& Topology Publications, Coventry, 2002.

\bibitem{HillmanLinnell92}
{\bf Hillman, J.~A. and Linnell, P.~A.}:
\newblock \emph{``Elementary amenable groups of finite {H}irsch length are
  locally-finite by virtually-solvable"},
\newblock {\em J. Austral. Math. Soc. Ser. A}, 52(2):237--241, 1992.

\bibitem{Linnell(1993)} {\bf Linnell, P.}:
\emph{``Division rings and group von Neumann algebras''},
Forum Math. 5, 561--576 (1993).

\bibitem{Lyndon-Schupp(1977)}
{\bf Lyndon, R.C. and Schupp, P. E.}:
\emph{``Combinatorial group theory''},
Ergebnisse der Mathematik und ihrer Grenzgebiete 89,
Springer (1977).

\bibitem{Lueck(2002)} {\bf L\"uck, W.}:
\emph{``$L^2$-Invariants: Theory and Applications to Geometry and 
$K$-Theory''}, Ergebnisse der Mathematik und ihrer Grenzgebiete 44,
Springer (2002).

\bibitem{Schick(2000)} {\bf Schick, T.}:
\emph{``Integrality of $L^2$-Betti numbers''},
Math. Ann. 317, 727--750 (2000).

\bibitem{Stallings(1968)} {\bf Stallings, J.R.}:
\emph{``On torsion-free groups with infinitely many ends''},
Annals of Math. 88, 312--334 (1968).
\end{thebibliography}
\end{document}